\documentclass[a4paper,reqno]{amsart}
\usepackage{amssymb, amsmath, amscd}

\newtheorem{theorem}{Theorem}[section]

\newtheorem{remark}[theorem]{Remark}
\newtheorem{definition}[theorem]{Definition}

\newtheorem{example}[theorem]{Example}
\def\spec{\operatorname{Spec}}
\def\Id{\operatorname{Id}}
\def\TJ{\tilde{\mathcal{J}}}
\def\Span{\operatorname{Span}}
\makeatletter
 \@addtoreset{equation}{section}
\makeatother
\begin{document}
\title{Projective Affine Ossermann Curvature Models}
\author{Peter Gilkey \text{and} Bronson Lim}
\address{PG: Mathematics Department, \; University of Oregon, \;\; Eugene \; OR
97403 \; USA}
\email{gilkey@uoregon.edu}
\address{BL: Department of Mathematics, University of Oregon, Eugene, OR 97403, USA}
\email{bcl@uoregon.edu}

\begin{abstract} A curvature model $ (V,\mathcal{A})$ is a vector space equipped with an element 
$A\in V^\ast\otimes V^\ast\otimes\text{End}(V)$ such that $A$ has the same symmetries as an 
affine curvature operator. 
Such a model is called projective affine Osserman if the spectrum of the Jacobi operator, 
$\mathcal{J}_X(y) = \mathcal{A}(y,x)x$, 
is projectively constant. There are topological conditions imposed by 
Adam's Theorem (vector fields on spheres) on such a model. In this paper 
we construct projective affine Osserman curvature models in dimensions 
$m\equiv1(2)$, $m\equiv 2(4)$, and $M\equiv4(8)$ which realize all possible 
eigenvalue structures allowed by Adam's Theorem.
\end{abstract}

\maketitle

\section{Introduction}
An {\it affine manifold} is a pair $\mathcal{A}:=(M,\nabla)$ where $M$ is a smooth 
$m$-dimensional manifold and where
$\nabla$ is a torsion free connection on the tangent bundle $TM$ of $M$.
Let $\mathcal{R}$ be the associated {\it curvature operator} which is defined by:
$$\mathcal{R}(X,Y):=[\nabla_X,\nabla_Y]-\nabla_{[X,Y]}\,.$$
The curvature operator has the following universal symmetries:
\begin{equation}\label{E1.a}
\mathcal{R}(X,Y)=-\mathcal{R}(Y,X)\text{ and }
\mathcal{R}(X,Y)Z+\mathcal{R}(Y,Z)X+\mathcal{R}(Z,X)Y=0\,.
\end{equation}
Note that any other universal symmetries of the curvature operator
in the affine setting are algebraic consequences of these \cite{GSW09}.

The {\it Jacobi operator}
$\mathcal{J}_X:Y\rightarrow\mathcal{R}(Y,X)X$
plays an important role in the study of geodesic sprays. We consider the
spectrum of $\mathcal{J}_X$:
$$\spec\{\mathcal{J}_X\}:=
\{\lambda\in\mathbb{C}:\det(\mathcal{J}_X-\lambda\Id)=0\}\,.
$$
It is natural to examine the interplay between the geometry of the manifold and the spectrum of
the Jacobi operator. We say
that $(M,\nabla)$ is {\it projective affine Osserman} if $\spec\{\mathcal{J}_X\}$
is a projective invariant which is independent of the particular non-zero tangent vector
$X$ which is chosen -- we refer to Definition~\ref{D1.3} for details.

We shall begin by putting matters in a historical perspective. 
In Section~\ref{S1.1}, we discuss the Osserman conjecture in
pseudo-Riemannian geometry. In Section~\ref{S1.2}, we present some material concerning
affine Osserman manifolds. In Section~\ref{S1.3}, we review previous results concerning
projective affine Osserman manifolds. 

In the original study of the Osserman conjecture in the
Riemannian setting, one first proceeded algebraically and used methods 
of algebraic topology to restrict the possible eigenvalue structures. 
One then exhibited purely algebraic Riemann curvature models using
Clifford modules to show that all the possible eigenvalue structures existed
algebraically.  Chi \cite{C88}  and Nikolayevskey \cite{niko1}--\cite{niko5} 
then used the second Bianchi identity to eliminate
many of the purely algebraic structures
and determine exactly which could be realized geometrically; the case 
$m=16$ had to be excluded at least in part since there was an Osserman 
Riemannian curvature model (given by the Cayley plane) not of Clifford type.
With this in mind,  in Section~\ref{S1.4}, 
we shall define the notion of a projective affine curvature model; this is a purely
algebraic object which encodes the relevant geometric condition at a single point. 
We will then state the main result of the paper (Theorem~\ref{T1.11}) 
which controls the possible eigenvalue structures if 
$m\equiv1(2)$, $m\equiv 2(4)$, and $m\equiv 4(8)$. 
The remainder of the paper is devoted to the proof of this result which
restricts the possible geometries which can arise and which we hope
will lead subsequently to a better understanding of projective affine 
Osserman geometry.

\subsection{The Osserman conjecture in pseudo-Riemannian Geometry}\label{S1.1}
Let $\nabla$ be the Levi-Civita connection of a pseudo-Riemannian manifold 
$\mathcal{M}=(M,g)$. Let
$$S^\pm(\mathcal{M}):=\{\xi\in TM:g(\xi,\xi)=\pm1\}$$
be the pseudo-sphere bundles of unit  timelike ($-$) and unit spacelike ($+$) 
tangent vectors.  Motivated by the seminal work of Osserman \cite{oss},
one says that $\mathcal{M}$ is 
{\it spacelike Osserman} (resp. {\it timelike Osserman})  if the spectrum of 
$\mathcal{J}_X$ is constant on $S^+(\mathcal{M}$) (resp. on $S^-(\mathcal{M})$). 
As these are equivalent conditions \cite{GKVa}, 
one simply says the manifold in question is {\it Osserman}.

Restrict to the Riemannian setting for the moment. 
One says that a simply connected complete
Riemannian manifold $\mathcal{M}$
is a $2$-point homogeneous space if the isometry group of $\mathcal{M}$ 
acts transitively on 
$S^+(\mathcal{M})$. Such a manifold is either flat or is a rank $1$-symmetric space,
i.e. the sphere, complex projective space,
quaternionic projective space, the Cayley plane, and the negative curvature duals.
For such a manifold, clearly $\spec\{\mathcal{J}_X\}$ is
independent of the particular $X\in S^+(\mathcal{M})$ which is chosen. 
Osserman \cite{oss} wondered if
the converse was true; 
if $\mathcal{M}$ is a simply connected complete Riemannian manifold and if 
$\spec(\mathcal{J}_X)$ is constant on $S^+(\mathcal{M})$, 
can one then conclude that 
$\mathcal{M}$ is a two point homogeneous space.
For this reason, his name has been associated with this area.  
This question has been answered in the affirmative by work of Chi \cite{C88}
and Nikolayevskey \cite{niko1}--\cite{niko5} except in dimension $m=16$;
the situation if $m=16$ still is not settled although there are some partial results.

A Lorentzian manifold is Osserman if and only if it has
constant sectional curvature; thus the classification is complete in this 
setting \cite{BBG97,GKVa}. 
In the higher signature setting, the situation is much more complicated. 
In dimension $4$, the only case to be considered is in signature $(2,2)$;
there is an intimate relationship between signature $(2,2)$ Osserman manifolds 
and Walker geometry
(see for example the discussion in \cite{D08} and associated references).
The algebraic structure is completely classified and there are 4 basic types. The
Jacobi operator can be diagonalizable, or can have a complex eigenvalue, or can
have real eigenvalues with a Jordan block of size $2$, or can
have real eigenvalues with a Jordan block of size 3.
At the geometric level, complex eigenvalues can not occur so this is an example
of an algebraic possibility which does not appear geometrically. There
are a great many results in this area and we can only cite a few in the interests of 
brevity \cite{CGV09,D09, DGV1}. In the
higher signature setting, the Jordan normal form can be arbitrarily complicated 
on the algebraic level \cite{GI01}. 
If $p\ge2$ and if $q\ge2$, then
there are Osserman pseudo-Riemannian manifolds of signature $(p,q)$
which are not locally symmetric\index{symmetric} \cite{GVV}.
The Rak\'ic duality principle \cite{AR07} has been studied in the 
pseudo-Riemannian setting  \cite{BM12,R99}. 
In the Riemannian setting, the quaterions
play an important role; para-quaternion (i.e. hypercomplex) structures play a
 corresponding role in neutral-signature Osserman geometry \cite{A05}.

Osserman geometry is intimately related with many questions in mathematical
physics. Chaichi et. al. \cite{CGM05} studied conditions for a Walker metric
to be Einstein, Osserman, or locally conformally flat and obtained thereby 
exact solutions to the Einstein equations for a
restricted Walker manifold. Chudecki and Prazanowski \cite{ChuPrz1,ChuPrz2} 
examined Osserman metrics in terms of $2$-spinors
and provided some new results in HH-geometry using the close relation between 
weak HH-spaces and Walker
and Osserman spaces using results of \cite{DGV06}.
One can also use the Weyl conformal curvature operator to define a conformal 
analogue of the Jacobi operator and
study the resulting geometry \cite{BG05,N12}. The geometry of the skew-symmetric 
curvature operator has been
studied analogously \cite{Cal09,C-GR-VL-07,IP2-98} as has the geometry of the
higher order Jacobi operator \cite{SV92}, of the higher order skew-symmetric
curvature operator
\cite{S04}, and of the Szab\'o operator \cite{Szabo}. In addition, other properties
of natural operators defined by the curvature can be examined \cite{BG07}; the field is a broad and fertile one.

\subsection{Affine Osserman manifolds}\label{S1.2}
In the affine setting, there is no metric to normalize the choice of tangent vector.
Since $\mathcal{J}_{\lambda X}=\lambda^2\mathcal{J}_X$, 
it is natural to say that an affine manifold $(M,\nabla)$ is
{\it affine Osserman} if $\mathcal{J}_X$ is nilpotent for all $X$ or, equivalently, 
if $\spec(\mathcal{J}_X)=\{0\}$
for all $X$. Such manifolds arise naturally
as generalized affine plane wave manifolds (see Section 2.2 \cite{GGNV13}).
This notion was first explored by E. Garc\'{\i}a-R\'{\i}o et. al. \cite{GKVV} and 
has proven to be a fruitful field
of inquiry. If $(M,\nabla)$ is
an affine manifold, let $(x^1,\dots,x^m)$ be local coordinates on $M$. 
If $\omega\in T^*M$, expand $\omega=\sum_iy_idx^i$ to define the dual fiber
coordinates $(y_1,\dots,y_m)$ and thereby obtain
{\it canonical local coordinates} $(x^1,\dots,x^m,y_1,\dots,y_m)$ 
on the cotangent bundle $T^*M$. 
Let $\Phi=\Phi_{ij}dx^i\circ dx^j$ be a
smooth symmetric $2$-tensor on $M$. 
The {\it deformed Riemannian extension} $g_{\nabla,\Phi}$
is the metric of neutral signature on the cotangent bundle $T^*M$ given locally by
\begin{equation}\label{E1.b}
\begin{array}{l}
g_{\nabla,\Phi}(\partial_{x_i},\partial_{x_j})=-2y_k\Gamma_{ij}{}^k(x)+\Phi_{ij}(x),\\
g_{\nabla,\Phi}(\partial_{x_i},\partial_{y^j})=\delta_i^j,\quad
g_{\nabla,\Phi}(\partial_{y^i},\partial_{y^j})=0\,.
\end{array}\end{equation}
This is invariantly defined -- see, for example, the discussion in \cite{CGGV09}. 
One has (see for example Theorem 2.16 of \cite{GGNV13}) that:
\begin{theorem}
Let $(M,\nabla)$ be an affine Osserman manifold and let $\Phi$ be a smooth 
symmetric $2$-tensor on $M$. Then $(T^*M,g_{\nabla,\Phi})$
is a pseudo-Riemannian nilpotent Osserman manifold of neutral signature.
\end{theorem}

The {\it modified Riemannian extension} is the neutral signature metric on the cotangent
bundle given locally by
$$\tilde g:= 2\, dx^i\circ dy_i+\{y_iy_j-2y_k\Gamma_{ij}{}^k\}dx^i\circ dx^j\,.$$
Again, this is
invariantly defined and one has \cite{CGGV12}:
\begin{theorem}
If $(M,\nabla)$ is an affine Osserman manifold, then $(T^*M,\tilde g)$ is
an Osserman manifold of neutral signature and the
eigenvalues of the Jacobi operator
on the unit pseudo-sphere\index{pseudo-sphere} bundles $S^
\pm(T^*M,\tilde g)$ are $\pm(0,1,\frac14)$ with multiplicities
$(1,1,2m-2)$, respectively.
\end{theorem}

If $(M,\nabla)$ is flat, then $(T^*M,\tilde g)$ is locally isomorphic to para-complex
projective space with constant para-holomorphic sectional curvature $+1$. 
Taking $(M,\nabla)$ to be non-flat gives rise to Osserman metrics on
neutral signature manifolds with non-nilpotent Jacobi operators and 
with non-trivial Jordan normal form which admit natural para-Hermitian structures. 
They are semi para-complex space forms which neither satisfy
the third Gray identity nor need they be integrable \cite{CGGV12}.

\subsection{Projective affine Osserman manifolds}\label{S1.3}
Another way to deal with the problem of rescaling in the affine category is to
projectivize the question.  Since $\mathcal{J}_XX=0$, zero is always in the spectrum
of $\mathcal{J}_X$. We assume that $(M,\nabla)$ is not affine Osserman and thus
that $\spec\{\mathcal{J}_X\}\ne\{0\}$ for some $X$.
We follow the discussion in \cite{BGNS08} and make the following:

\begin{definition}\label{D1.3}
\rm\ \begin{enumerate}
\item  An affine manifold $ (M,\nabla) $ is said to be {\it projective affine
Osserman} if $(M,\nabla)$ is not affine Osserman and if
there is a smooth non-vanishing function $s$ so 
$\spec\{\mathcal{J}_X\} =
s(Y,X)\spec\{\mathcal{J}_Y\}$ for all $X,Y\in TM-\{0\}$.
\item A pseudo-Riemannian manifold $(M,g)$ is said to be {\it projective
spacelike} ($+$) or {\it timelike} ($-$) {\it Osserman}
if $(M,g)$ is not flat and if
there exists a smooth non-vanishing function $s_\pm$ so that 
$\spec\{\mathcal{J}_X\}
=s_\pm(Y,X)\spec\{\mathcal{J}_Y\}$ for all $X,Y$ in $S^\pm(M,g)$.
\end{enumerate}
\end{definition}

Adopt the notation of Equation~(\ref{E1.b}) to define the
deformed Riemannian extension $(T^*M,g_{\nabla,\Phi})$. We have \cite{BGNS08}:

\begin{theorem}\label{T1.4}
\ \begin{enumerate}
\item If $(M,g)$ is a Riemannian Osserman manifold which is not flat
and if $\nabla$ is the Levi-Civita connection, then $(M,\nabla)$
is projective affine Osserman.
\item If$(M,\nabla)$ is an affine manifold,
then following assertions are equivalent:
\begin{enumerate}
\item $(M,\nabla)$ is a projective affine Osserman manifold.
\item $(T^*M,g_{\nabla,\Phi})$ is a projective spacelike Osserman manifold.
\item $(T^*M,g_{\nabla,\Phi})$ is a projective timelike Osserman manifold.
\end{enumerate}\end{enumerate}
\end{theorem}

\begin{remark}\rm
We note that  indefinite signature Osserman metrics are not projective affine
Osserman as  $\spec\{\mathcal{J}(X)\}=0$ if $X$ is null.
\end{remark}

To illustrate this, we present two
examples of projective affine Osserman manifolds that
do not arise from an underlying Riemannian structure and refer
to \cite{gs1} for other examples:

\begin{example}\label{E1.6}
\rm Let $\{e_1,...,e_m\}$ be the standard basis for $\mathbb{R}^m$.
Expand $x=x^ie_i$ to define the usual coordinates $(x^1,\dots,x^{m})$ 
on $\mathbb{R}^{m}$. 
We suppose $m\ge2$ and let $\mathcal{M}:=(\mathbb{R}^m,\nabla)$ where the
Christoffel symbols of $\nabla$ are given by:
$$
\Gamma_{mm}{}^m=2,\ \Gamma_{im}{}^i=\Gamma_{mi}{}^i=\Gamma_{ii}{}^m=1\text{ for }i<m;
\ \Gamma_{11}{}^1=-\Gamma_{22}{}^2=\varepsilon(x_1+x_2)
$$
where $\epsilon\in\mathbb{R}$ is a real parameter.
We have
$$R_{ijk}{}^l=\partial_{x_i}\Gamma_{jk}{}^l-\partial_j\Gamma_{ik}{}^l+\Gamma_{in}{}^l\Gamma_{jk}{}^n-
\Gamma_{jn}{}^l\Gamma_{ik}{}^n\,.$$
There are no terms in $\varepsilon^2$ and the only terms in $\varepsilon$ which are quadratic
in the Christoffel symbols are
\begin{eqnarray*}
0&=&\Gamma_{m1}{}^1\Gamma_{11}{}^1-\Gamma_{11}{}^1\Gamma_{m1}{}^1=0,\\
0&=&\Gamma_{m2}{}^2\Gamma_{22}{}^2-\Gamma_{22}{}^2\Gamma_{m2}{}^2=0\,.
\end{eqnarray*}
Consequently, the quadratic terms give rise to:
\begin{eqnarray*}
&&R_{imm}{}^i=\Gamma_{im}{}^i\Gamma_{mm}{}^m-\Gamma_{mi}{}^i\Gamma_{im}{}^i
=2-1\text{ for }i<m,\\
&&R_{mii}{}^m=\Gamma_{mm}{}^m\Gamma_{ii}{}^m-\Gamma_{ii}{}^m\Gamma_{mi}{}^i
=2-1\text{ for }i<m,\\
&&R_{ijj}{}^i=\Gamma_{im}{}^i\Gamma_{jj}{}^m=1\text{ for }i\ne j<m\,.
\end{eqnarray*}
We also have
$$R_{122}{}^2=\partial_{x_1}\Gamma_{22}{}^2=-\varepsilon,\quad
R_{211}{}^1=\partial_{x_2}\Gamma_{11}{}^1=\varepsilon\,.$$
One may then verify \cite{gs1} that this is projective affine Osserman. The
entries in the curvature tensor are constant so this manifold is affine curvature
homogeneous.

Suppose $\epsilon\ne0$. Then
$\mathcal{J}(e_3)$ is diagonal. If
$X=\frac1{\sqrt2}(e_1+e_3)$, then:
\begin{eqnarray*}
&&\mathcal{J}(x)(e_1+e_3)=0,\quad
\mathcal{J}(x)(e_1-e_3)=e_1-e_3,\quad
\mathcal{J}(x)e_2=\textstyle\frac12\varepsilon e_1+e_2,\\
&&\mathcal{J}(x)\{e_2+\textstyle\frac14\varepsilon(e_1+e_3)\}=\frac12\varepsilon e_1+e_2
=e_2+\textstyle\frac14\varepsilon(e_1+e_3)+\textstyle\frac14\varepsilon(e_1-e_3)\,.
\end{eqnarray*}
Thus $\Span\{v_1:=e_2+\textstyle\frac14\varepsilon(e_1+e_3),v_2=(e_1-e_3)\}$
is invariant under the action of $\mathcal{J}(x)$. As $\mathcal{J}(x)v_2=v_2$ and
$\mathcal{J}(x)v_1=v_1+v_2$, we have non-trivial Jordan normal form in this instance.

There is a translation group of rank $m-1$ which acts on $(M,\nabla)$
preserving the structures. We have additional entries in $\nabla R$:
$$
\nabla R(\partial_2,\partial_1,\partial_1;\partial_1)
=-2\Gamma_{11}{}^1\partial_2,\text{ and }
\nabla R(\partial_1,\partial_2,\partial_2;\partial_2)=-2\Gamma_{22}{}^2\partial_1\,.
$$
Since $\Gamma_{11}{}^1$ and $\Gamma_{22}{}^2$ vanish if and only if $x_1+x_2=0$,
$(M,\nabla)$ is not 1-affine curvature homogeneous and has
affine cohomogeneity 1.

Suppose $\varepsilon=0$.
The curvature tensor is then that of constant sectional curvature.
Since the Christoffel symbols are constant, the group of translations acts
transitively on $\mathcal{M}$ by affine isomorphisms; thus $\mathcal{M}$ is affine
homogeneous\index{homogeneous}. However, if we set $\sigma(t)=(0,...,0,x(t))$, then the
geodesic equation\index{geodesic equation}
becomes $\ddot x+2\dot x\dot x=0$
which blows up in finite time for suitable initial conditions. Thus
$(\mathbb{R}^m,\nabla)$ is geodesically incomplete\index{geodesically incomplete}. Finally, 
we compute:
\begin{eqnarray*}
&&\nabla R(\partial_m,\partial_1,\partial_1;\partial_m)\\
&=&\nabla_{\partial_m}R(\partial_m,\partial_1)\partial_1
-R(\nabla_{\partial_m}\partial_m,\partial_1)\partial_1-R(\partial_m,\nabla_{\partial_m}\partial_1)\partial_1
-R(\partial_m,\partial_1)\nabla_{\partial_m}\partial_1\\
&=&(2-2-2)\partial_m\ne0\,.
\end{eqnarray*}
Consequently, $\nabla R\ne0$. Thus this manifold
is not locally symmetric\index{symmetric} and
is not affinely equivalent to the standard affine structure on the sphere $S^m$.
\end{example}

\begin{example}\label{E1.7}
\rm Let $(M,\nabla)$ be an affine surface. Let $\rho_s$ be symmetric part of the Ricci tensor
defined by $\nabla$. Because $\mathcal{J}_XX=\mathcal{R}(X,X)X=0$ and
because the dimension is two, one of the following two possibilities pertains:
\begin{enumerate}
\item $\spec\{\mathcal{J}_X\}=\{0\}$ and $\rho(X,X)=0$.
\item $\spec\{\mathcal{J}_X\}=\{0,\lambda(X)\}$ for $0\ne\lambda(X)\in\mathbb{R}$
 and $\rho(X,X)=\lambda(X)\ne0$.
 \end{enumerate}
Thus  $\rho$ is definite if and only if $\rho(X,X)\ne0$ for $X\ne0$
or, equivalently, if $\spec\{\mathcal{J}_X\}\ne\{0\}$ for $X\ne0$.
Thus $(M,\nabla)$ is projective affine Osserman if and only $\rho_s$ is definite. Since the
alternating Ricci tensor $\rho_a$ can be non-trivial, this provides many examples of projective
affine Osserman surfaces where the connection is not the Levi-Civita connection.
\end{example}

\subsection{Curvature Models}\label{S1.4} It is convenient to work purely in the
algebraic setting to understand what happens at each tangent space. 
The corresponding analysis in the Riemannian setting was central to the classification
in dimensions $m\ne16$. 

\begin{definition}
\rm Let $V$ be a finite dimensional real vector space of dimension $m$.
If $\mathcal{A}\in V^\ast\otimes V^\ast\otimes\text{End}(V)$ has the curvature symmetries of
Equation~(\ref{E1.a}), then the pair $(V,\mathcal{A})$ is said to be an {\it affine curvature model} and 
the associated {\it Jacobi operator} is given by $\mathcal{J}_X(Y):= \mathcal{A}(Y,X)X$.
An affine manifold $(M,\nabla)$ is said to be a
{\it geometric realization of an affine curvature model}
$(V,\mathcal{A})$ if there is a point $P$ of $M$ and a linear isomorphism $\psi:T_PM\rightarrow V$ so that 
$\phi^*\mathcal{A}=\mathcal{R}_P$. 
\end{definition}
Every affine curvature model is geometrically realizable \cite{GSW09}; 
thus, as noted above, the symmetries
of Equation~(\ref{E1.a}) generate the universal symmetries of the curvature operator of an affine manifold.

\begin{definition}\rm
We say that $ (V,\mathcal{A})$ is a {\it projective affine Osserman curvature model}
if $(V,\mathcal{A})$ is an affine curvature model and
if the spectrum of the Jacobi operator is a projective invariant, i.e.
if for every $0\ne X,Y\in V$ there exists $0\neq s(Y,X)\in\mathbb{R}$ so
$$\spec\{\mathcal{J}_X\} =s(Y,X)\spec\{\mathcal{J}_Y\}\neq \{0\}\,$$
If $\lambda\in\mathbb{C}$, let $\mu(\mathcal{J}_X,\lambda)=\mu(\lambda)$ be
the eigenvalue multiplicity; $\lambda\in\spec\{\mathcal{J}_X\}$ if and only
if $\mu(\mathcal{J}_X,\lambda)>0$. As $\mathcal{J}_X$ is real, $\mu(\lambda)=\mu(\bar\lambda)$.
\end{definition}

In fact, the nature of the spectrum can be nailed down a bit if $(V,\mathcal{A})$
is a projective affine Osserman curvature model. 
In the following result, we shall omit the $\lambda$'s 
if there are no real eigenvalues and the $\mu$'s if there are no complex eigenvalues.
We refer to \cite{gs1} for proof of:
\begin{theorem} \label{T1.10}
Let $(V,\mathcal{A})$ be a projective affine Osserman curvature model. 
There exist smooth complex valued non-zero functions
$\{\lambda_1(X),\dots,\mu_1(X),\dots\}$ defined on $V-\{0\}$ taking distinct values
and there exists a smooth positive function $s(\cdot,\cdot)$ defined for $0\ne X,Y\in V$ so that:
\begin{enumerate}
\item $\lambda_i(X)\in\mathbb{R}-\{0\}$ and 
$\mu_j(X)\in\mathbb{C}$ with $\Im(\mu_j)>0$.
\item  $\spec\{\mathcal{J}_X\}
=\{0,\lambda_1(X),\dots,\nu_1(X),\bar\nu_1(X),\dots\}$.
\item $\lambda_i(X)=s(X,Y)\lambda_i(Y)$ and $\nu_j(X)=s(X,Y)\nu_j(Y)$
 for all $i$ and $j$.
\item Let
$\vec\mu(X):=(\mu(0),\mu(\lambda_1),\dots,\mu(\nu_1),\mu(\bar\nu_1),\dots)$
where $\mu(\cdot)$ is the eigenvalue multiplicity. Then $\vec\mu(\cdot)$ is constant
on $V-\{0\}$.
\end{enumerate}
\end{theorem}

Since $\mathcal{J}_XX=0$, $0$ is always in the spectrum of the Jacobi operator.
It is convenient therefore to introduce the {\it reduced Jacobi operator}. Let
$V_X:=V/\{X\cdot\mathbb{R}\}$. Since $X\in\ker\{\mathcal{J}_X\}$, the Jacobi
operator induces a natural map $\TJ_X$, which is called the
{\it reduced Jacobi operator}, from $V_X$ to $V_X$. 
Let $\mu(\mathcal{J}_X,\sigma)$ and $\mu(\TJ_X,\sigma)$ be the
eigenvalue multiplicities of $\sigma$ in $\mathcal{J}_X$ or $\TJ_X$. As
$\det(\mathcal{J}_X-\lambda\Id)=-\lambda\det(\TJ_X-\lambda\Id)$,
\begin{eqnarray*}
&&\spec(\mathcal{J}_X)=\spec(\TJ_X)\cup\{0\}\text{ and }\\
&&\mu(\mathcal{J}_X,\sigma)=\left\{\begin{array}{lll}
\mu(\TJ_X,0)+1&\text{if}&\sigma=0\\
\mu(\TJ_X,\sigma)&\text{if}&\sigma\ne0\end{array}\right\}\,.
\end{eqnarray*}

The following is the main result of this paper and deals with the cases 
$m\equiv1(2)$, $m\equiv2(4)$ and $m\equiv4(8)$. 
Note that we permit $\lambda_i(X)=0$ for some $i$ in what follows except in (1), (2a), or (3a).

\begin{theorem}\label{T1.11}
Let $ (V,\mathcal{A})$ be a projective affine Ossermann curvature model of dimension $m$. Let $X\ne0$.
\begin{enumerate}
\item Suppose $m\equiv1(2)$. Then $\spec\{\TJ_X\}=\{\lambda_1(X)\}$ and
$\vec\mu=(m-1)$.
\smallbreak\item Suppose $m\equiv2(4)$. Then one of the following possibilities holds:
\begin{enumerate}
\item $\spec\{\TJ_X\}=\{\lambda_1(X)\}$ and
$\vec\mu=(m-1)$.
\item $\spec\{\TJ_X\}=\{\lambda_1(X),\lambda_2(X)\}$ and
$\vec\mu=(1,m-2)$.
\item $\spec\{\TJ_X\}=\{\lambda_1(X),\nu_1(X),\bar\nu_1(X)\}$ and
$\vec\nu=(1,\frac{m-2}2,\frac{m-2}2)$.
\end{enumerate}
\smallbreak\item Suppose $m\equiv4(8)$.  Then one of the following possibilities holds:
\begin{enumerate}
\item $\spec\{\TJ_X\}=\{\lambda_1(X)\}$ and 
$\vec\mu=(m-1)$.
\item $\spec\{\TJ_X\}=\{\lambda_1(X),\lambda_2(X)\}$ and
{\rm[i]} $\vec\mu=(1,m-2)$ or
{\rm[ii]} $(2,m-3)$\par\qquad or {\rm[iii]} $\vec\mu=(3,m-4)\}$.
\item $\spec\{\TJ_X\}=\{\lambda_1(X),\lambda_2(X),\lambda_3(X)\}$ and
{\rm[i]} $\vec\mu=(1,1,m-3)$,\par\qquad or {\rm[ii]} $\vec\mu=(1,2,m-4)\}$.
\item $\spec\{\TJ_X\}=\{\lambda_1(X),\lambda_2(X),\lambda_3(X),\lambda_4(X)\}$
and $\vec\mu=(1,1,1,m-4)$.
\item $\spec\{\TJ_X\}=\{\lambda_1(X),\nu_1(X),\bar\nu_1(X)\}$ and
{\rm [i]} $\vec\mu=(1,\frac{m-2}2,\frac{m-2}2),$\par\qquad or {\rm[ii]}
$\vec\mu=(3,\frac{m-4}2,\frac{m-4}2)$, or {\rm[iii]} $\vec\mu=(m-3,1,1)$.
 \item $\spec\{\TJ_X\}=\{\lambda_1(X),\lambda_2(X),\nu(X),\bar{\nu}(X)\}$
 and {\rm[i]} $\vec\nu=(1,2,\frac{m-4}2,\frac{m-4}2)$,\par\qquad
 or {\rm[ii]} $(1,m-4,1,1)\}$.
\item $\spec\{\TJ_X\}
=\{\lambda_1(X),\lambda_2(X),\lambda_3(X),\nu(X),\bar{\nu}(X)\}$ and
\par\qquad
$\vec\mu=(1,1,1,\frac{m-4}2,\frac{m-4}2)$.
\item $\spec\{\TJ_X=(\lambda_1(X),\nu_1(X),\bar\nu_1(X),\nu_2(X),\bar\nu_2(X))$ and\par\qquad
$\vec\mu=(1,1,1,\frac{m-4}2,\frac{m-4}2)$.
\end{enumerate}
\item There is a projective affine Osserman curvature model realizing
each of the eigenvalue structures above.
\end{enumerate}\end{theorem}

If $m=2$, then only case (2-a) appears.
Similarly, if $m=4$, then many of the possibilities are not present.

In Section \ref{S2}, we recall Adam's Theorem on vector fields on spheres and
use it to restrict the possible eigenvalue structures to establish Assertions~(1)--(3).
In Section~\ref{S3}, we will prove Assertion~(4) and
show that Assertions~(1)--(3) are sharp by constructing
examples which realize all the indicated structures.

\section{Methods of algebraic topology}\label{S2}

\subsection{Adam's Theorem} Chi \cite{C88} noticed that one could attack the
Osserman problem by looking at decompositions of the tangent bundle of the
sphere, $ TS^{m-1}$, into sub-bundles. Since decompositions of
$ TS^{m-1}$ correspond to linearly independent vector fields on the
sphere, the following is an immediate consequence of work of Adams \cite{A62}:

\begin{theorem}\label{T2.1}
Let $S^{m-1}$ be the unit sphere in $\mathbb{R}^m$. Suppose we can decompose
the tangent bundle $TS^m=E_1\oplus E_2\oplus...\oplus E_\ell$ as the direct sum of
vector bundles where $1\le\dim(E_1)\le...\le\dim(E_\ell)$.
\begin{enumerate}
\item If $m\equiv1(2)$, then $\ell=1$ and $\dim(E_\ell)=m-1$.
\item If $m\equiv2(4)$, then $\ell\le 2$ and $\dim(E_\ell)\ge m-2$.
\item If $m\equiv4(8)$, then $\ell\le 4$ and $\dim(E_\ell)\ge m-4$.
\end{enumerate}
\end{theorem}

This result is sharp. Let $E_0(X):=X\cdot\mathbb{R}$.  
If $m\equiv2(4)$, then we may regard $\mathbb{R}^m=\mathbb{C}^{m/2}$. Let
$E_1(X):=\sqrt{-1}X\cdot\mathbb{R}$ and $E_2(X):=(E_0(X)\oplus E_1(X))^\perp$.
This gives a decomposition of $T(S^{m-1})=E_0(X)^\perp=E_1\oplus E_2$ where
$\dim(E_2)=m-2$. If $m\equiv4(8)$, then we may regard
$\mathbb{R}^m=\mathbb{H}^{m/4}$ where $\mathbb{H}=1\cdot\mathbb{R}\oplus
I\cdot\mathbb{R}\oplus J\cdot\mathbb{R}\oplus K\cdot\mathbb{R}$ 
is the skew-field of quaternions. Setting 
$E_1(X):=IX\cdot\mathbb{R}$, $E_2(X):=JX\cdot\mathbb{R}$, 
$E_3(X):=KX\cdot\mathbb{R}$, and
$E_4(X):=(E_0(X)\oplus E_1(X)\oplus E_2(X)\oplus E_3(X))^{\perp}$ 
then produces a decomposition
with $\ell=4$ and $\dim(E_\ell)=m-4$.

\subsection{The proof of Theorem~\ref{T1.11}~(1,2,3)} Let $(V,\mathcal{A})$ be a
projective affine Osserman curvature model with associated
Jacobi operator $\mathcal{J}$. Let $V_X:=V/\{X\cdot\mathbb{R}\}$
and let $\TJ_X:V_X\rightarrow V_X$ be the
{\it reduced Jacobi operator} defined previously. We will work with 
$\TJ$ to prove Theorem~\ref{T1.11}.
Endow $V$ with an auxiliary positive definite inner product
$\langle\cdot,\cdot\rangle$. Let
$$S =S(V):=\{X\in V: \langle X,X\rangle = 1\}$$ be the associated sphere.
Since vector spaces $T_XS$ and $V/\{\mathbb{R}\cdot X\}$ may be canonically identified,
the reduced Jacobi operator $\TJ_X$ defines an endomorphism of the
tangent bundle $TS$.
If $\sigma$ is a (possibly) complex eigenvalue, let
$$
E_\sigma(X):=\{\xi\in V/\{\mathbb{R}\cdot X\}:
(\TJ_X-\sigma)^m(\TJ_X-\bar\sigma)^m\xi=0\}
$$
be the generalized eigenspaces of $\TJ_X$. 
We use Theorem~\ref{T1.10} to enumerate
the eigenspaces $E_{\lambda_1},\dots$ and $E_{\nu_1},\dots$; 
since the multiplicities are constant, these
patch together to define smooth vector bundles which define a decomposition
$$T(S^{m-1})=E_{\lambda_1}\oplus\dots\oplus E_{\nu_1}\oplus\dots\text{ where }\Im(\nu_i)>0\,.$$

If $m\equiv1(2)$, then there is one eigenbundle. 
One uses an argument using characteristic
classes (see \cite{gs1}) to rule out the case that the eigenbundle relates to a 
complex eigenvalue
and conclude that the eigenvalue $\lambda$ is real; then an example may be
obtained by taking a space of constant sectional curvature $\lambda$ or, if 
positive, by using a
rescaled version of the manifold given in Example~\ref{E1.6}. 
This completes the proof of Theorem~\ref{T1.11}
if $m\equiv1(2)$.

If $m\equiv2(4)$ and if there is only one bundle in the decomposition of 
$T(S^{m-1})$, then the corresponding eigenvalue must be real as 
$\dim(T(S^{m-1}))$ is odd; this gives rise to Case (2-a). 
If there are two bundles in the decomposition, then one must be a line bundle 
and thus corresponds to a real eigenvalue $\lambda_1$. 
Case (2=b) arises when the complementary bundle corresponds to a real eigenvalue 
$\lambda_2$ and Case (2-c) arises when the complementary bundle corresponds to a complex eigenbundle.
The analysis of the case $m\equiv4(8)$ is similar and is therefore omitted. 
This completes the proof of
the first 3 assertions of Theorem~\ref{T1.11}.\hfill\qed

\section{Projective affine Osserman curvature models}\label{S3}
We complete the proof of Theorem~\ref{T1.11} by constructing projective
affine Osserman curvature models which realize each of the eigenvalue
structures given above. Let $\langle\cdot,\cdot\rangle$ be the usual inner product
on $\mathbb{R}^m$. We shall be considering two different basic operators. 
Let $\mathcal{A}_0$ be the curvature operator of constant sectional curvature:
$$
\mathcal{A}_0(X,Y)Z:=\langle Y,Z\rangle X-\langle X,Z\rangle Y\,.
$$
We verify that $\mathcal{A}_0$ is an affine algebraic curvature operator by computing:
\medbreak\qquad
$\mathcal{A}_0(X,Y)Z=\langle Y,Z\rangle X-\langle X,Z\rangle Y=-
\mathcal{A}_0(Y,X)Z$,
\smallbreak\qquad
$\mathcal{A}_0(X,Y)Z+\mathcal{A}_0(Y,Z)X+\mathcal{A}_0(Z,X)Y$
\smallbreak\qquad\quad
$=\langle Y,Z\rangle X-\langle X,Z\rangle Y$
$+\langle Z,X\rangle Y-\langle Y,X\rangle Z$
$+\langle X,Y\rangle Z-\langle Z,Y\rangle X$\smallbreak\quad\qquad$=0$.
\medbreak\noindent
Next suppose that $m$ is even and that $J$ is a unitary almost complex structure, i.e.
$J^2=-\Id$ and $J^*\langle\cdot,\cdot\rangle=\langle\cdot,\cdot\rangle$. We define
$$\mathcal{A}_J(X,Y)Z:={\textstyle\frac13}\{\langle JY,Z\rangle X-\langle JX,Z\rangle Y-2\langle JX,Y\rangle Z\}\,.$$
We use the relation $\langle JX,Y\rangle=-\langle X,JY\rangle$ to show that
$\mathcal{A}_J$ is an affine curvature operator by computing:
\medbreak\quad\qquad
$\mathcal{A}_J(X,Y)Z:=
\frac13\{\langle JY,Z\rangle X-\langle JX,Z\rangle Y-2\langle JX,Y\rangle Z\}=-
\mathcal{A}_J(Y,X)Z$,
\smallbreak\quad\qquad
$\mathcal{A}_J(X,Y)Z+\mathcal{A}_J(Y,Z)X+\mathcal{A}_J(Z,X)Y$
\smallbreak\qquad\qquad\qquad
$=\frac13\{\langle JY,Z\rangle X-\langle JX,Z\rangle Y-2\langle JX,Y\rangle Z\}$
\smallbreak\qquad\qquad\qquad
$+\frac13\{\langle JZ,X\rangle Y-\langle JY,X\rangle Z-2\langle JY,Z\rangle X\}$
\smallbreak\qquad\qquad\qquad
$+\frac13\{\langle JX,Y\rangle Z-\langle JZ,Y\rangle X-2\langle JZ,X\rangle Y\}=0$.
\medbreak\noindent
If $\Xi$ is an auxiliary linear operator, then $\Xi\mathcal{A}_0$ and $\Xi\mathcal{A}_J$ are affine
algebraic curvature operators. If $\|X\|=1$, then the Jacobi operators are given by:
\begin{equation}\label{E3.a}
\mathcal{J}^{\Xi\mathcal{A}_0}_XY=\left\{\begin{array}{l}
\Xi Y\text{ if }Y\perp X\\
\quad0\text{ if }X=Y
\end{array}\right\}\text{ and }
\mathcal{J}_X^{\Xi\mathcal{A}_J}Y=\left\{
\begin{array}{l}\Xi X\text{ if }Y=JX\\
\quad0\text{ if }Y\perp JX\end{array}\right\}\,.
\end{equation}
We note that $\mathcal{A}_0+\frac13J\mathcal{A}_J$ is the operator
of constant holomorphic sectional curvature.
\subsection*{Step I: $m$ arbitrary} 

Let $\mathcal{A}:=a_0\mathcal{A}_0$. If $\|X\|=1$, then
$\spec\{\Xi J_X\}=\{a_0\}$; this constructs a projective affine Osserman
curvature model realizing the eigenvalue structures given in (1), (2-a), and (3-a) of Theorem~\ref{T1.11}.

\subsection*{Step II: $m\equiv0\mod 2$} 
Let 
\begin{equation}\label{E3.b}
\textstyle\mathcal{A}=a_0\mathcal{A}_0+a_1J(\mathcal{A}_0-J\mathcal{A}_J)+(c_1-a_0)J\mathcal{A}_{\mathcal J}\,.
\end{equation}
Let $\|X\|^2=1$. The associated Jacobi operator is given by
$$\mathcal{J}_X^{\mathcal{A}}Y=\left\{\begin{array}{rl}
0&\text{if }Y=X\\
c_1Y&\text{if }Y=JX\\
a_0Y+a_1JY&\text{if }Y\perp X,JX
\end{array}\right\}\,.$$
Let $\lambda_1\in\mathbb{R}$ and $\nu_1\in\mathbb{C}-\mathbb{R}$ with $\Im(\mu_1)>0$. Let
$(c_1,a_0,a_1)=(\lambda_1,\Re(\nu_1),\Im(\nu_1))$. Then
$$
\spec\{\\XiJ_X^{\mathcal{A}}\}=\{\lambda_1,\mu_1,\bar\mu_1\}\text{ and }
\vec\mu=(1,(m-2)/2,(m-2)/2)
$$
 which is case (2-c). If $\lambda_1\ne\lambda_2$,
 take $(c_1,a_0,a_1)=(\lambda_1,\lambda_2,0)$. Then
 $$
 \spec\{\\XiJ_X\}=\{\lambda_1,\lambda_2\}\text{ and }\vec\mu=(1,m-2)$$
 which is
case (2-b) and (3-b-i). Taking $\lambda_1=\lambda_2$ yields case (2-a).

\subsection*{Step III: $m\equiv0\mod4$} 
Identify $\mathbb{R}^m=\mathbb{H}^{m/4}$ to define a quaternion
structure $\{J_1,J_2,J_3\}$ on $\mathbb{R}^m$ where $J_1J_2=J_3$.
Let $\{\lambda_1,\lambda_2,\lambda_3,\lambda_4\}$ be real. Set
\begin{equation}\label{E3.c}
\begin{array}{l}
\mathcal{A}:=c_4\mathcal{A}_0+(c_1-a_0)J_1\mathcal{A}_{J_1}
+(c_2-a_0)J_2\mathcal{A}_{J_2}+(c_3-a_0)J_3\mathcal{A}_{J_3}\\
\qquad+a_1J_1(\mathcal{A}_0-J_1\mathcal{A}_{J_1})+a_2J_1(J_2\mathcal{J}_2+J_3\mathcal{J}_3)\,.
\vphantom{\vrule height 11pt}\end{array}\end{equation}
We use Equation~(\ref{E3.a}) to see
$$\mathcal{J}_X^{\mathcal{A}}Y=\left\{\begin{array}{lll}
0&\text{ if }Y=X\\
c_1Y&\text{ if }Y=J_1X\\ 
c_2Y+(a_1+a_2)J_1Y&\text{ if }Y=J_2X\\ 
c_3Y+(a_1+a_2)J_1Y&\text{ if }Y=J_3X\\ 
c_4Y+a_1J_1Y&\text{ if }Y\perp\{X,J_1X,J_2X,J_3X\}\end{array}\right\}\,.
$$
Since $J_1J_2J_3=\Id$, $J_1$ defines a complex structure which preserves the spaces
$\Span\{J_2X,J_3X\}$ and $\Span\{X,J_1X,J_2X,J_3X\}^\perp$; this is an essential point that fails
for higher rank Clifford module structures and which prevents us extending this construction to the
case $m\equiv8\mod16$. We consider the following cases:
\begin{enumerate}
\item Assume all the eigenvalues are real. Let $c_i=\lambda_i$, $a_1=0$, and $a_2=0$. Then 
$\mathcal{J}_X^{\mathcal{A}}$ is diagonalizable. If all the eigenvalues are distinct, we obtain the
structure of (3-d); the structures of (3-a), (3-b), and (3-c) are obtained by letting the eigenvalues coalesce.
\item Assume there is one complex eigenvalue $\nu_1$ of multiplicity $(m-4)/2$. We let $c_4=\Re(\nu_1)$,
$a_1=\Im(\nu_1)$, $a_2=-a_1$, $c_1=\lambda_1$, $c_2=\lambda_2$, and $c_3=\lambda_3$ to obtain
the eigenvalue structure in (3-g); the eigenvalue structures in (3-e-ii) and (3-f-i) arise by letting the
real eigenvalues coalesce.
\item Assume there is one complex eigenvalue $\nu_1$ of multiplicity $1$ and two real eigenvalues 
$\lambda_1$ and $\lambda_2$. We let $a_1=0$, $c_1=\lambda_1$,
$c_2=c_3=\Re(\mu_1)$, $c_4=\lambda_2$, and $a_2=\Im(\mu_1)$ to obtain the eigenvalue structure of (3-f-ii); the
eigenvalue structure of (3-e-iii) arises from taking $\lambda_1=\lambda_2$.
\item Assume there is one complex eigenvalue $\nu_1$ of multiplicity $(m-2)/2$ and one real eigenvalue $\lambda_1$
of multiplicity $1$. We take $\lambda_1=c_1$, $c_2=c_3=c_4=\Re(\mu_1)$, $a_1=\Im(\mu_1)$, and $a_2=0$
to obtain the eigenvalue structure of (3-e-i).
\item Assume there are two distinct complex eigenvalues $\nu_1$ and $\nu_2$ and one real eigenvalue
$\lambda_1$. Take $c_1=\lambda_1$, $c_2=c_3=\Re(\nu_1)$, $a_1+a_2=\Im(\nu_1)$, 
$c_4=\Re(\nu_2)$, and $a_1=\Im(\nu_2)$
to obtain the eigenvalue structure of (3-h).
\end{enumerate}
This completes the proof of Theorem~\ref{T1.11}.\hfill\qed

\begin{remark}\rm
 There are Lie groups underlying the constructions we have given
to prove Assertion~(4) of Theorem~\ref{T1.11}. Let $\operatorname{GL}(\mathbb{R}^m)$
be the general linear group of invertible linear maps from $\mathbb{R}^m$
to $\mathbb{R}^m$. The {\it orthogonal group} is defined by setting:
$$
\mathcal{O}:=\{\Xi\in\operatorname{GL}(\mathbb{R}^m):\Xi^*\langle\cdot,\cdot\rangle=
\langle\cdot,\cdot\rangle\}\,.
$$
The operator $\mathcal{A}_0$ of constant sectional curvature is invariant under
the action of $\mathcal{O}$ since it arises out of the quadratic form
$\langle\cdot,\cdot\rangle$. Since $\mathcal{O}$ acts transitively on $S^{m-1}$,
$\mathcal{A}_0$ is a projective affine Oserman operator. If $m$ is even,
let $J$ be a Hermitian almost complex structure on
$\mathbb{R}^m=\mathbb{C}^{m/2}$. The associated {\it unitary group} is given by:
$$\mathcal{U}:=\{\Xi\in\mathcal{O}:J\Xi=\Xi J\}\,.$$
The curvature operator given in Equation~(\ref{E3.b}) is invariant under the action of $\mathcal{U}$.
Since $\mathcal{U}$ acts transitively on $S^{m-1}$,
the curvature operator of Equation~(\ref{E3.a})
is a projective affine Osserman curvature operator. Finally, suppose $m$ is
divisible by $4$. We identify $\mathbb{R}^m=\mathbb{H}^{m/4}$ to give
$\mathbb{R}^m$ a quaternion structure. The associated {\it quaternion} group
is given by:
$$\mathcal{SP}:=\{\Xi\in\mathcal{O}:J_1\Xi=\Xi J_1,\text{ and }J_2\Xi=\Xi J_2\}\,.$$
This group acts transitively on $S^{m-1}$ and preserves the
curvature operator of Equation~(\ref{E3.c}); since $\mathcal{SP}$ acts transitively on
$S^{m-1}$, the curvature operator of Equation~(\ref{E3.c}) is projective affine Osserman.
Unfortunately, this process terminates at this stage; 
higher order Clifford module structures
do not give rise to transitive group actions on spheres.
\end{remark}

\subsection*{Acknowledgements}
Research partially supported by project MTM2009-07756 (Spain).

\end{document}